\documentclass[12pt]{article}
\usepackage{amsmath}
\usepackage{amsfonts}
\usepackage{amssymb}
\usepackage{amsthm}
\usepackage[english]{babel}
\usepackage{epic,eepic}
\usepackage{fancyhdr}
\usepackage{indentfirst}

\makeatletter
\@addtoreset{equation}{section}
\makeatother

\newtheorem{statement}{}[section]
\newtheorem{definition}[statement]{Definition}
\newtheorem{theorem}[statement]{Theorem}
\newtheorem{proposition}[statement]{Proposition}
\newtheorem{corol}[statement]{Corollary}
\newtheorem{lemma}[statement]{Lemma}

\def\noi{\noindent}
\newcommand{\C}{\mathbb C}

\newcommand{\Z}{\mathbb Z}
\newcommand{\N}{\mathbb N}
\newcommand{\T}{\mathbb T}
\newcommand{\E}{\mathbb E}
\newcommand\e{{\rm e}}
\newcommand{\eps}{\varepsilon}
\newcommand\ind{{\rm 1\kern-.30em I}}

\newcommand{\biindice}[3]%
{%

\begin{array}[t]{c}
{\displaystyle #1}\\
{\scriptstyle #2}\\
{\scriptstyle #3}
\end{array}

}

\begin{document}

\title{\bf Thin sets of integers in Harmonic analysis and $p$-stable random Fourier series}
\author{\it P. Lef\`evre, D. Li, H. Queff\'elec, L. Rodr{\'\i}guez-Piazza\\
\\
{\normalsize \it Dedicated to the memory of Daniel Rider}
}
\date{ \footnotesize \today}

\maketitle

\bigskip

\noi {\bf Abstract.} {\it We investigate the behavior of some thin sets of integers defined through random trigonometric 
polynomial when one replaces Gaussian or Rademacher variables by $p$-stable ones, with $1 < p< 2$. We show 
that in one case this behavior is essentially the same as in the Gaussian case, whereas in another case, this 
behavior is entirely different.}
\medskip

\noi {\bf AMS Classification numbers.} Primary: 43A46 ; secondary: 42A55; 42A61; 60G52 
\medskip

\noi{\bf Key-words.} Fourier series; quasi-independent set; Rider set; stable random variables; stationary set 

\section{Introduction} 

Let $G$ be a compact, abelian group (which will be mostly the circle $\T$), equipped with its normalized Haar 
measure $m$, and $\Gamma$ its (discrete) dual. We will denote by  ${\cal P}$ the set of finite sums 
$\sum_{\gamma\in\Gamma} c_{\gamma}\gamma$, i.e. the vector space generated by $\Gamma$, and by 
${\cal P}_{\Lambda}$ the set of finite sums $\sum_{\gamma\in\Lambda} c_{\gamma}\gamma$, where 
$\Lambda$ is a subset of $\Gamma$. We recall (\cite{Ru}) that $\Lambda$ is called a Sidon set if, for some 
constant $C$, we have the following \emph{a priori} inequality:
\begin{equation} \label{1.1}
\qquad\qquad
\Vert f\Vert_{F_1} := \sum_\gamma\vert \hat f(\gamma)\vert\leq C\Vert f\Vert_\infty, 
\quad \forall f\in {\cal P}_\Lambda.
\end{equation} 

The best constant $C$ in \eqref{1.1} is called the Sidon constant of $\Lambda$.
A long standing problem, solved in the positive  by Drury (\cite{Dru}) at the beginning of the seventies, was 
whether the union of two Sidon sets is again a Sidon set. A little after Drury's result, Rider (\cite{Ri}) gave 
the following necessary and sufficient condition for  Sidonicity, from which the result becomes obvious:
\begin{equation}\label{1.2}
\Vert f\Vert_{F_1}\leq C\, [\![f]\!], \quad \forall f\in {\cal P}_\Lambda.
\end{equation}
Here, we have set:
\begin{equation}\label{1.3} 
[\![f]\!] = {\mathbb E}\, \Big\Vert \sum \varepsilon_{\gamma} \hat f (\gamma) \gamma \Big\Vert_\infty,
\end{equation}  
where $(\varepsilon_\gamma)_\gamma$ is a sequence of i.i.d. Rademacher random variables, defined 
on some probability space  $\Omega$, i.e. independent and taking the values $+1$ and $-1$ with equal 
probability $\frac{1}{2}$\,, and where $\E$ stands for expectation on $\Omega$. This norm was thoroughly 
studied by Marcus and Pisier  (\cite{MarPis})  and is called the ${\cal C}^{as}$-norm in the space of almost 
surely continuous random Fourier series. These two authors proved in particular the non-trivial fact that one 
could as well use a standard gaussian sequence instead of a Rademacher one, and obtain an equivalent norm 
(see \cite{Pi}, Th\'eor\`eme~7.1). Pisier (\cite{Pi}) realized that Sidonicity can also  be characterized by the 
\emph{a priori} inequality: 
\begin{equation}\label{1.4}
\Vert f\Vert _{\infty}\leq C\, [\![f]\!], \quad \forall f \in {\cal P}_\Lambda.
\end{equation}

This is a general fact, the proof of which we recall for the convenience of the reader, and which motivates the 
forthcoming definition of  \emph{stationarity}: let $(Z_\gamma)_{\gamma\in\Gamma}$ be a collection of 
i.i.d. copies of a complex-valued, centered and integrable random variable $Z$, and set, for every 
trigonometric polynomial $f$, 
\begin{displaymath}
[\![ f ]\!]_{_Z} = \E\, \Big\Vert \sum Z_{\gamma}\hat f (\gamma) \gamma \Big\Vert_\infty.
\end{displaymath}
We have the following simple proposition.

\begin{proposition} \label{prop Sidon}
Let $\Lambda\subset \Gamma$ be such that 
\begin{equation}\label{Sidon}
\qquad 
\Vert f \Vert_\infty \leq C\,  [\![ f ]\!]_{_Z} \,, \quad \forall f\in{\cal P}_{\Lambda}.
\end{equation}
Then, $ \Lambda$ is a Sidon set. 
\end{proposition}
 
\noi{\bf Proof.} Let  $(\widetilde Z_\gamma)_{\gamma\in\Gamma})$ be an independent family, with each 
$\widetilde Z_{\gamma}$ a symmetrization of $Z_\gamma$. Since the latter variables are centered, we have   
\begin{displaymath}
\E\, \Big\Vert \sum Z_{\gamma}\hat f (\gamma) \gamma \Big\Vert_\infty 
\leq 2\, \E\,\Big\Vert \sum \widetilde Z_{\gamma} \hat f (\gamma) \gamma \Big\Vert_\infty,
\end{displaymath}
so we may as well assume the $Z_\gamma$'s symmetric from the beginning. If 
$\varepsilon_\gamma= \pm 1$, we therefore have
\begin{displaymath}
 \Big\vert \sum_{\gamma} \varepsilon_\gamma \hat f (\gamma) \Big\vert 
\leq \Big\Vert \sum_{\gamma} \varepsilon_\gamma \hat f (\gamma)\gamma \Big\Vert_\infty
\leq C\, \E\, \Big\Vert \sum Z_{\gamma}\hat f (\gamma) \gamma \Big\Vert_\infty.
\end{displaymath}
Taking the supremum on all choices of $\pm 1$ gives classically (\cite{LiQ}, Chapitre~5, Proposition~IV.2):
\begin{displaymath}
\sum_{\gamma}\vert\hat f (\gamma)\vert 
\leq D\, \E\, \Big\Vert \sum Z_{\gamma}\hat f (\gamma) \gamma \Big\Vert_\infty, 
\end{displaymath}
with $D = \frac{\pi}{2}\,C$. Now, we truncate our variables $Z_\gamma$ at a level $M$:  
$Z_\gamma = Z'_\gamma + Z''_\gamma$, where 
\begin{displaymath}
Z'_\gamma = Z_\gamma \ind_{\{\vert Z_\gamma\vert \leq M \}} ,
\quad  Z''_\gamma = Z_\gamma \ind_{\{\vert Z_\gamma\vert > M \}}; 
\end{displaymath}
$M$ being adjusted so as to have $\E\,\vert Z''_\gamma\vert\leq \frac{1}{2D}$\,. We now see that 
\begin{align*}
\sum_\gamma \vert\hat f (\gamma)\vert 
&\leq D\, \E\, \Big\Vert \sum_\gamma Z'_{\gamma}\hat f (\gamma) \gamma \Big\Vert_\infty 
+ D\, \E\, \Big\Vert \sum_{\gamma} Z''_{\gamma}\hat f (\gamma) \gamma \Big\Vert_\infty  \\
&  \leq D\, \E\, \Big\Vert \sum_\gamma Z'_{\gamma}\hat f (\gamma) \gamma \Big\Vert_\infty 
+ D\, \sum_\gamma \vert \hat f (\gamma)\vert \,\E\, \vert Z''_\gamma\vert \\
& \leq D\, \E\, \Big\Vert \sum_\gamma Z'_{\gamma} \hat f (\gamma)  \gamma \Big\Vert_\infty 
+ \frac{1}{2} \sum_\gamma\vert \hat f (\gamma)\vert \\
& \leq 4MD\, \E\, \Big\Vert\sum_\gamma \varepsilon_\gamma \hat f (\gamma)\gamma \Big\Vert_\infty  
+ \frac{1}{2} \sum_\gamma\vert \hat f (\gamma)\vert
\end{align*}
(here, $(\varepsilon_\gamma)_\gamma$ is a Rademacher sequence, and we used the usual 
``contraction principle": see \cite{LiQ}, Chapitre~3, Th\'eor\`eme~III.3); whence 
\begin{displaymath}
\sum_\gamma\vert \hat f (\gamma)\vert 
\leq 8MD\, \E\, \Big\Vert\sum_{\gamma}\varepsilon_\gamma \hat f (\gamma)\gamma \Big\Vert_\infty. 
\end{displaymath}

Now, we are in position to apply Rider's Theorem (\cite{LiQ}, Chapitre~5, Th\'eor\`eme~IV.18) to conclude 
that $\Lambda$ is Sidon. 
\hfill $\square$
\medskip

Note that Proposition~\ref{prop Sidon} has an easy converse (which we state for further reference).
\begin{proposition}\label{converse Sidon}
If $\Lambda$ is a Sidon set, then $\| \,\, \|_\infty$ and $[\![\,\,]\!]_{_Z}$ are equivalent norms on 
${\cal P}_\Lambda$.
\end{proposition}
{\bf Proof.} Let  $f \in {\cal P}_\Lambda$ and 
$f_{_Z}^\omega = \sum_\gamma Z_\gamma (\omega) \widehat f (\gamma)\,\gamma$. On one hand, one has
\begin{displaymath}
\E_\omega \| f\|_\infty \leq \sum_{\gamma \in \Lambda} \E (|Z_\gamma|)\,|\widehat f (\gamma) |
= \|Z_1\|_1 \sum_{\gamma \in \Lambda} |\widehat f (\gamma) |\leq C\,\|Z \|_1 \|f \|_\infty,
\end{displaymath}
and, on the other hand, 
$\sum_\gamma | Z_\gamma (\omega) \widehat f (\gamma) | \leq C\, \|f_{_Z}^\omega\|_\infty$, 
since $\Lambda$ is Sidon and $f_{_Z}^\omega \in {\cal P}_\Lambda$; hence, by integrating,  
$\sum_\gamma \E (| Z_\gamma|) \,|\widehat f (\gamma)| \leq C\,[\![ f]\!]_{_Z}$, and
\begin{equation}
\| f \|_\infty \leq \sum_\gamma |\widehat f (\gamma) | \leq \frac{C}{\| Z_1\|_1} \, [\![ f ]\!]_{_Z}.
\tag*{$\square$}
\end{equation}
\medskip

Pisier (\cite{Pi})  also studied the subsets $\Lambda$ of $\Gamma$ verifying the reverse inequality of 
\eqref{Sidon} (see \cite{Pi}, D\'efinition~6.2), namely:
\begin{equation}\label{1.5}
 [\![f]\!]\leq C\, \Vert f \Vert _{\infty} , \quad \forall f\in {\cal P}_\Lambda.
\end{equation}
and he called those sets  \emph{stationary}, proving in particular (\cite{Pi}, Proposition~6.2) that the 
cartesian product of $d$ Sidon sets is always stationary (the first named author \cite{Le} proved that, for 
example,  $\{3^{k_1} + \cdots + 3^{k_d}\,;\ 1 \leq k_1 < \cdots < k_d\}$ is also a stationary set). Another 
well-known notion is that of $q$-Sidonicity, $1\leq q < 2$ (\cite{LoRo}, \cite{Pi}, D\'efinition~6.1). 
The subset $\Lambda$ is called \emph{$q$-Sidon} if, for some constant $C$, we have:
\begin{equation}\label{1.6}
\qquad
\Vert f\Vert_{F_q} :=\Big(\sum_{\gamma}\vert\hat f (\gamma)\vert^q\Big)^{\frac{1}{q}} 
\leq C\,\Vert f\Vert_\infty, \quad \forall f\in{\cal P}_{\Lambda}.
\end{equation}

After the work of Rider, the following notion was also introduced (\cite{LiQuRo} and \cite{Ro}): the subset 
$\Lambda$ is called \emph{$q$-Rider} if, for some constant $C$, we have this time:
\begin{equation}\label{1.7}
\qquad 
\Vert f\Vert_{F_q}\leq C\, [\![f]\!], \qquad \forall f\in{\cal P}_\Lambda.
\end{equation}

It is immediate to see that every $q$-Sidon set is a $q$-Rider set, and Rider's result can be formulated in 
saying that the converse holds for $q=1$. Whether this converse holds for each $q\in (1,2)$ is an open 
problem, in spite of several non-trivial partial results (\cite{LeRo}). Let us mention that the cartesian 
product of $d$ infinite Sidon sets is $q$-Sidon with $q = \frac{2d}{d+1}$ and not better (\cite{LoRo}). \par
\smallskip

We could of course study those notions for other probability laws than the (subgaussian) Rademacher laws 
or gaussian ones. This is precisely the aim of this work, where we will be interested in the complex, symmetric, 
$p$-stable random variables $Z$, $1 < p < 2$, which can be defined through their characteristic function:
\begin{displaymath}
\qquad\qquad
\E\,\big(\e^{i{\rm R}e\,\overline{z}Z}\big) = \exp(-\vert z\vert^p),\quad \forall z\in\C.
\end{displaymath}

The case $p=2$ is the gaussian case already studied. The case $1 < p < 2$ is in some sense more delicate, 
because in spite of the nice \emph{stability} property: 
\begin{displaymath}
\sum_{n=1}^N a_n Z_n \sim \Big(\sum_{n=1}^N\vert a_n\vert^p \Big)^{\frac{1}{p}} Z_1 
\end{displaymath}
from which those variables borrow their name, their integrability properties  are fairly poor: $Z \in L^s$ 
for each $s < p$, but $Z\not\in L^p$ (in fact, $Z\in L^{p,\infty}$). Yet, this case  has also  been studied in 
great detail by Marcus and Pisier in \cite{MaPi}, who in particular introduced the following $p$-stable  norm 
on the space of trigonometric polynomials
\begin{equation}\label{1.8}
\qquad \qquad
[\![f]\!]_p = \E\, \Big\Vert \sum_{\gamma\in\gamma} 
Z_\gamma \hat f (\gamma)\gamma \Big\Vert_\infty,\quad \forall f\in{\cal P},
\end{equation}
where $(Z_\gamma)_\gamma$ is a family of independent copies of a complex  $p$-stable, symmetric, 
random variables. Observe that this has a meaning,  since the $Z_\gamma$'s are integrable. Moreover, due 
to a general comparison principle of Jain and Marcus (\cite{JaMa}), one has the following inequality, where 
the implied constants only depend on $p_1$ and $p_2$:
\begin{equation}\label{1.9}
\qquad \quad
1< p_1 < p_2\leq 2 \quad\Rightarrow \quad [\![f]\!]_{p_2}\leq C (p_1,p_2)\, [\![f]\!]_{p_1},
\ \forall   f\in{\cal P}.
\end{equation}

In other terms, the smaller $p$, the bigger the corresponding $[\![ \ ]\!]_{p}$-norm. In particular, those new 
norms are bigger than the previously mentioned Rademacher and Gaussian norms on ${\cal P}$. The 
questions which we examine in this work are the following: what do the notions of stationarity, 
$q$-Riderness, become if we replace the gaussian variables by $p$-stable ones? \par

\bigskip

After having established, in Section~2, one basic property of the $[\![ \  ]\!]_{p}$-norm, namely a lower 
$p$-estimate, we prove in Section~3, that a $p$-stationary set is in fact Sidon (and, of course, conversely) 
as soon as $p < 2$, and we  study, in Section~4, several equivalent forms of $p$-stable $q$-Riderness, and 
show that this apparently new notion coincides with that of $s$-Riderness for an appropriate value of the 
parameter $s$, depending on $p$ and $q$. We end with some comments.\par


\section{ Basic properties of the $p$-stable norm}

We will need the following two theorems on  $p$-stable norms. They are more or less straightforward 
consequences of a basic result of Marcus and Pisier. \par
First, we  introduce some notation: $F_p$ will denote the set of functions $f \in L^2= L^2 (G,m)$ such that 
their Fourier transform is in $\ell_p = \ell_p (\Gamma)$, equipped with the norm 
\begin{displaymath}
\Vert f\Vert_{F_p} :=\Vert \hat f\Vert_p \,, 
\end{displaymath}
which we already encountered in Section 1 (see \eqref{1.6}). 
\smallskip

We shall denote by $\|.\|_\psi$ the Luxemburg norm in the Orlicz space associated to an Orlicz function $\psi$. Let $r>0$, we shall be mainly interested in the Orlicz function $\varphi_r$, where
\begin{displaymath}
\varphi_{r} (x) = x\big (1+\log (1+x) \big)^{\frac{1}{r}},
\end{displaymath}
and the conjugate Orlicz function $\psi_r$, where
\begin{displaymath}
\psi_{r} (x) = {\rm e}^{x^r}-1.
\end{displaymath}
\smallskip

Finally let $A (p,\varphi_{p'})$ (where $p'$ is the conjugate exponent of 
$p$) be the space of all functions in $L^2(G)$ which can be written as 
\begin{displaymath}
f = \sum_{n=1}^\infty h_n\ast k_n \,,
\end{displaymath}
with:
\begin{displaymath}
\sum_{n=1}^\infty\Vert h_n\Vert_{F_p} \Vert k_n\Vert_{\varphi_{p'}} <\infty,
\end{displaymath} 
and equipped with the norm 
\begin{displaymath}
\Vert f\Vert_{A(p,\varphi_{p'})} = 
\inf\Big\{\sum_{n=1}^\infty \Vert h_n\Vert_{F_p}\Vert k_n\Vert_{\varphi_{p'}}\Big\}\,,
\end{displaymath}
where the infimum runs over all possible representations of $f$.\par
\smallskip\goodbreak

With those notations, the  basic result alluded to above stands as follows, under a simplified form which 
will be sufficient for us (\cite{MaPi}, Theorem~5.1):

\begin{theorem}[{\bf Marcus-Pisier}] \label{theo Marcus-Pisier}
The norms  $[\![f]\!]_{p}$ and $\Vert f\Vert_{A(p,\varphi_{p'})}$ are equivalent on the space ${\cal P}$ of 
trigonometric polynomials on $G$.
\end{theorem}

Two important consequences of that theorem, which are not explicited in \cite{MaPi} in the $p$-stable case, 
are the following (see \cite{Pi}, Proposition~7.1 for the case $p = 2$).

\begin{theorem}[{\bf Contraction principle for the $p$-stable norm}] 
Let $f\in {\cal P}$ and $(\xi_\gamma)_\gamma$ be a collection of functions in $L^p(0,1)$, bounded in 
$L^p$. Denote the integral over $(0,1)$ by $\E'$. Then, we have for some  positive constant $a$:
\begin{equation}\label{2.1}
\Big(\E'\, \Big[\!\!\Big[\sum_{\gamma}\xi_\gamma 
\hat f (\gamma)\gamma \Big]\!\!\Big]_p^p\Big)^{\frac{1}{p}}
\leq a \sup_{\gamma} \Vert \xi_\gamma\Vert_{L^p(0,1)} \, [\![f]\!]_{p}.
\end{equation} 
\end{theorem}
\noi{\bf Proof.} Let  $f =\sum_{n=1}^\infty h_n\ast k_n$  be an admissible decomposition of $f$, and 
let $\omega'\in (0,1)$, as well as 
$f_{\omega'} =\sum_{\gamma} \xi_{\gamma} (\omega') \hat f (\gamma) \gamma $. We can write 
$f_{\omega'} =\sum_{n=1}^\infty H_{n} (\omega') \ast k_n$, with $H_{n} (\omega')\in F_p$ and 
$\widehat {H_{n} (\omega')} (\gamma) =\xi_{\gamma} (\omega')\hat h_{n}(\gamma)$. Set for convenience 
$X (\omega') = [\![ f_{\omega'} ]\!]_{p}$ and $Y_{n} (\omega') =\Vert H_n (\omega') \Vert_{F_p}$.
We see that 
\begin{displaymath}
X (\omega')\leq \sum_{n=1}^\infty [\![ H_{n}(\omega')\ast k_n ]\!]_{p} 
\leq a \sum_{n=1}^\infty Y_{n} (\omega') \Vert k_n\Vert_{\varphi_{p'}} \,,
\end{displaymath}
where $a$ is some  constant given by the Marcus-Pisier Theorem above. Now taking $L^p$-norms in 
$L^p(0,1)$ and using the triangle inequality, we get:
\begin{displaymath}
\Vert X \Vert_{L^p(0,1)}\leq a\sum_{n=1}^\infty \Vert Y_n \Vert_{L^p(0,1)} \Vert k_n\Vert_{\varphi_{p'}} 
\leq a C \sum_{n=1}^\infty \Vert h_n\Vert_{F_p} \Vert k_n\Vert_{\varphi_{p'}} \,,
\end{displaymath}
where $C = \sup_{\gamma}\Vert\xi_{\gamma}\Vert_{L^p(0,1)}$. Taking the infimum over all 
possible representations of $f$ gives us the result, possibly changing the constant $a$. \hfill $\square$

\bigskip

The second basic consequence is: 
\goodbreak

\begin{theorem} [\hbox{\bf Lower $p$-estimate for the $[\![ \,\  ]\!]_{p}$-norm}] \label{lower p-estimate}
 Let $f, f_1, \ldots, f_N $ be trigonometric polynomials such that:
\begin{displaymath}
\qquad \qquad
\vert \hat f (\gamma)\vert 
\geq \Big(\sum_{j=1}^N\vert \hat f_j (\gamma)\vert ^p \Big)^{\frac{1}{p}}, 
\qquad\forall \gamma\in\Gamma. 
\end{displaymath} 
Then, the constant $a$ being as in \eqref{2.1}:
\begin{equation}\label{2.2} 
[\![ f ]\!]_{p} \geq a^{-1} \Big(\sum_{j=1}^N[\![ f_j ]\!]_{p}^p \Big)^{\frac{1}{p}}\cdot
\end{equation}
\end{theorem}

\noi {\bf Proof.} Let $A_1, \ldots, A_N$ be a partition of $(0,1)$ in sets of Lebesgue measure $1/N$ and,  
for each $\gamma\in \Gamma$, $\xi_{\gamma}\in L^{p}(0,1)$ be defined by:  
\begin{displaymath}
\xi_{\gamma} = N^{\frac{1}{p}} \sum_{j=1}^N \frac{\hat f_j (\gamma)}{\hat f (\gamma)} \, \ind_{A_j}. 
\end{displaymath}

It is clear that $\Vert \xi_\gamma\Vert_{L^p(0,1)}\leq 1$ by our assumption, and by definition we have 

\begin{displaymath}
\E'\, \Big[\!\!\Big[ \sum_{\gamma}\xi_\gamma \hat f (\gamma)\gamma \Big]\!\!\Big]_p^p
=\sum_{j=1}^N  \Big[\!\!\Big[ \sum_{\gamma} 
\frac{\hat f_j (\gamma)}{\hat f (\gamma)}\hat f (\gamma) \gamma \Big]\!\!\Big]_{p}^p 
=\sum_{j=1}^N  [\![ f_j]\!]_{p}^p \,,
\end{displaymath}

\noi so that an application of \eqref{2.1} gives the result.
\hfill $\square$
\bigskip

We will end this section with the following estimate. This is undoubtedly known, but we did 
not find any explicit mention; so we are going to give some words of explanation. 

\begin{lemma}\label{majoriz. p-stable norm}
Let $\lambda_1<\lambda_2<\cdots<\lambda_n\in\N$, with $\lambda_n\ge2$. If $f (t) = \sum_{j=1}^n \e^{i {\lambda_j}t}$, $t\in \T$, one has, for some constant $C > 0$:
\begin{displaymath}
[\![ f ]\!]_p \leq C\,n^{1/p} (\log \lambda_n)^{1/p'}.
\end{displaymath} 
\end{lemma}
\noi{\bf Proof.} From (4.6) of \cite{MaPi} (or Remark~1.7, page~186 of \cite{MaPi2}), 
there is a constant $K > 0$ such that 
\begin{displaymath}
[\![ g ]\!]_p \leq K \sum_{k=2}^\infty \frac{1}{k (\log k)^{1/p}}\,
\Big( \sum_{j=k}^\infty |\widehat g (j) |^p \Big)^{1/p}
\end{displaymath}
for every trigonometric polynomial $g$ with spectrum in $\N$. Here one has 
$\sum_{j=k}^\infty |\widehat f (j) |^p \leq n$ for $k \leq \lambda_n$ and 
$\sum_{j=k}^\infty |\widehat f (j) |^p = 0$  for $k > \lambda_n$; hence
\begin{equation}
[\![ f ]\!]_p \leq 
K\,n^{1/p} \sum_{k=2}^{\lambda_n} \frac{1}{k (\log k)^{1/p}} \leq C\, n^{1/p} (\log \lambda_n)^{1/p'}.
\tag*{$\square$}
\end{equation}

\bigskip\goodbreak


\section{$p$-stable stationary sets are Sidon for $p < 2$}

The aim of this Section is to prove the following:

\begin{theorem} \label{theo p-stable stat = Sidon}
Let $\Lambda\subset \Gamma$ be a \emph{$p$-stationary} set ($1 < p \leq 2$), i.e. a set satisfying the 
following inequality, for some constant $C = C_p$:
\begin{displaymath}
\qquad
[\![ f ]\!]_p \leq C\, \Vert f \Vert_\infty, \quad \forall f\in {\cal P}_\Lambda.
\end{displaymath}
Then, if  $p < 2$, $\Lambda$ is a Sidon set.
\end{theorem}

The difference between the cases $1 < p <2$ and $p = 2$ (for which non-Sidon stationary sets exist), comes 
from the fact the $p$-stable norm for $p <2$ is bigger than the usual Pisier norm; hence having an upper 
estimate for it on some space forces the smallness of this space.\par
\bigskip

In order to prove Theorem~\ref{theo p-stable stat = Sidon}, we will need the following simple lemma :

\begin{lemma} \label{small lemma}
If $\Lambda\subset\Gamma$ is a $p$-stationary set, the norms  $[\![ \  ]\!]_2$  and  $[\![ \  ]\!]_p$ are 
equivalent on the space ${\cal P}_\Lambda$.
\end{lemma}

\noi {\bf Proof of Lemma~\ref{small lemma}.}
Let $f\in {\cal P}_\Lambda$, $(\varepsilon_\gamma)_\gamma$ a Rademacher sequence, and  
$f^{\omega} = \sum_\gamma \varepsilon_\gamma (\omega) \hat f (\gamma) \gamma$. By symmetry, we 
have $[\![ f  ]\!]_p = [\![ f^\omega  ]\!]_p  \leq C\, \Vert f^\omega\Vert_\infty$, where $C$ is the stationarity 
constant of $\Lambda$. Integrating over $\omega$ gives: 
$ [\![ f  ]\!]_p \leq C\,\E\, \Vert f^\omega\Vert_\infty = C [\![ f ]\!]_2$, which finishes the proof, 
since we know from \eqref{1.9} that the reverse inequality always holds (we wrote $ [\![ \ ]\!]_2$ instead of  
$[\![\  ]\!]$ to make a clear distinction between the subgaussian (i.e. Rademacher) or gaussian case, and the 
$p$-stable one).
\hfill $\square$

\medskip

\noi {\bf Proof of  Theorem~\ref{theo p-stable stat = Sidon}.} From the previous Lemma~\ref{small lemma}, 
it will be enough to show that, if $\Lambda$ is not Sidon, then the norms  $[\![ \ ]\!]_2 $ \ and $ [\![ \ ]\!]_p$ 
are not equivalent on ${\cal P}_\Lambda$. It will be convenient to introduce first some notation.  
Recall that a subset $B$ of $\Gamma$ is called \emph{quasi-independent} if, for any finite subset 
$\{\gamma_1, \ldots,\gamma_r\}$ of distinct elements of $B$, a relation 
$\sum_{i=1}^r \theta_i \gamma_i = 0$ , with $\theta_i = 0,\pm 1$ implies $\theta_i = 0$ for each $i$. 
The quasi-independent sets are the prototypes of Sidon sets, in that their Sidon constant $S$ is bounded by 
an absolute constant (\cite{LiQ}, Chapitre~12, Proposition~I.1, or \cite{LoRo}; in \cite{Kahane}, Kahane found 
that $S \leq 4.7$). Now, if $A$ is a finite subset of $\Gamma$:

\begin{enumerate} 
\item $\vert A\vert$ will denote the (finite) cardinality of $A$;
\item $q(A)$ will denote the biggest possible cardinality of a quasi-independent subset $B\subset A$;
\item $[\![A]\!]_2$  (respectively $[\![A]\!]_p$) will denote the quantity 
$\big[\!\!\big[ \sum_{\gamma\in A}\gamma \big]\!\!\big]_2$ (respectively  
$\big[\!\!\big[ \sum_{\gamma\in A}\gamma \big]\!\!\big]_p$);
\item $i_{A,p}$ will be the canonical injection of $({\cal P}_A,\Vert \ \Vert_{F_p})$ in $({\cal P}_A, [\![\ ]\!]_p$), 
and $\Vert i_{A,p}\Vert$ its norm.
\end{enumerate}

We will make use of the two following theorems, the first of which is an improvement of Rider's one, since it 
claims that it is enough to test the assumptions of that theorem on polynomials with coefficients $0$ or $1$.

\begin{theorem}[{\bf Pisier} \cite{Pis2}] 
The subset $\Lambda$ is Sidon if and only if there exists a constant $c>0$ such that
\begin{equation}\label{3.1} 
[\![A]\!]_2 \geq c\, \vert A\vert, \quad  {\rm for\ all\ finite }\ A\subset \Lambda.
\end{equation}
\end{theorem}

\noi (see \cite{LiQ}, Chapitre~12, Th\'eor\`eme~I.2).
\medskip\goodbreak

\begin{theorem}\label{extractgauss} ({\bf Rodr{\'\i}guez-Piazza} \cite{Ro2}; \cite{Ro}, Teorema~IV.1.3) 
There exists a numerical constant $K$ such that the following inequality holds :
\begin{equation}\label{3.2} 
K^{-1} q(A)\leq \Vert i_{A,2}\Vert^2\leq K q(A).
\end{equation} 
\noi In particular, we have :
\begin{equation}\label{3.3}
q(A) \geq K^{-1} \frac{[\![A]\!]_2^2}{\vert A\vert}\,\cdot
\end{equation}
\end{theorem}

\noi (see also \cite{LiQ}, Chapitre~12, Exercice~12.1, for a proof). To get \eqref{3.3} from \eqref{3.2}, it 
suffices to observe that $\Vert \sum_{\gamma\in A}\gamma\Vert _2^2=\vert A\vert$.
\medskip

Let us go back to the proof of  Theorem~\ref{theo p-stable stat = Sidon}. If $\Lambda$ is not Sidon, by 
\eqref{3.1}, we have:
\begin{displaymath}
\inf_{A\subset \Lambda,\,\vert A \vert <\infty} \frac{[\![A]\!]_2}{\vert A\vert} = 0, 
\end{displaymath} 
or equivalently: 
\begin{displaymath}
C_N :=\biindice{\inf}{A\subset\Lambda}{1\le\vert A\vert\le N} 
\frac{[\![A]\!]_2}{\vert A\vert} \mathop{\longrightarrow}_{N\to\infty}  0 .
\end{displaymath}

\medskip

Let now $\delta > 0$ and very small, and $N_0$ be the smallest integer such that $C_{N_0} <\delta$. Then, 
there exists a finite $A_0\subset \Lambda$, with $\vert A_0\vert\le N_0$ (and actually $\vert A_0\vert= N_0$ by definition of $N_0$) such that 
$\frac{[\![A_0]\!]_2}{\vert A_0\vert} <\delta$. We claim that we can find an integer $N$ and 
disjoint, quasi-independent subsets $B_1, \ldots, B_N$ of $A_0$ such that ($K$ being as in \eqref{3.2}) :
\begin{equation}\label{3.4}
\qquad \qquad\qquad\qquad
K^{-1} \delta\, [\![A_0]\!]_2 
\geq \vert B_j \vert \geq K^{-1} \frac{\delta}{2}\, [\![A_0]\!]_2, \quad  \forall \ 1\leq j \leq N;
\end{equation}
\begin{equation}\label{3.5}
\vert B_1\cup \cdots \cup B_N\vert \geq \frac{\vert A_0\vert}{2}\cdot
\end{equation}

Let us first see how \eqref{3.4} and \eqref{3.5} allow us  to finish the proof. They imply together:
\begin{displaymath}
\vert A_0\vert \leq 2\sum_{j=1}^N \vert B_j\vert \leq 2K^{-1} N\delta\, [\![A_0]\!]_2
\leq 2 K^{-1} N \delta^{2}\vert A_0\vert, 
\end{displaymath}
and so:
\begin{equation}\label{3.6} 
N \geq \frac{K}{2} \, \delta^{-2}.
\end{equation}

Now, the  lower $p$-estimate of Theorem~\ref{lower p-estimate} as well as the fact that 
$\vert B\vert \leq c_p^{-1} [\![B]\!]_p$ for quasi-independent sets (which are uniformy Sidon as we already 
mentioned), where $c_p$ is a  constant, gives us, $a$ being as in \eqref{2.2}:
\begin{align*}
[\![A_0]\!]_p^p 
& \geq a^{-p} \sum_{j=1}^N [\![B_j]\!]_p^p 
\geq a^{-p} c_p^p\sum_{j=1}^N \vert B_j \vert^p 
\geq a^{-p} \Big( \frac{c_p}{2K} \Big)^p N\delta^p [\![A_0]\!]_2^p \\
& \geq b_p \, \delta^{p-2} [\![A_0]\!]_2^p, 
\end{align*}
where $b_p$ only depends on $p$, and where we used \eqref{3.4} and \eqref{3.6}. Now, since $p < 2$ and since 
$\delta$ is arbitrarily small, this inequality proves that $[\![A_0]\!]_p$ is much bigger than  
$[\![A_0]\!]_2$  and ends the proof. \par
\medskip

It remains to show \eqref{3.4} and \eqref{3.5}. We first observe that 
\begin{equation}\label{3.7}
A' \subsetneqq A_0 \quad \Longrightarrow \quad [\![A']\!]_2\geq \delta \vert A'\vert , 
\end{equation}
and then:
\begin{equation}\label{3.7bis}
\frac{[\![A_0]\!]_2}{\vert A_0\vert}\geq \frac{\delta}{ 2} \,\cdot
\end{equation}
Indeed, the inequality \eqref{3.7} follows from the fact that $N_0$ is the smallest integer such that 
$C_{N_0} < \delta$. For the second inequality \eqref{3.7bis}, let $\emptyset \not= A' \subsetneqq A_0$, and set 
$A'' = A_0 \setminus A'$. The first inequality, applied to $A'$ and $A''$, as well as the unconditionality of the 
$[\![ \  ]\!]_2$-norm, give:
\begin{displaymath}
[\![A_0]\!]_2\geq [\![A']\!]_2\geq \delta \vert A'\vert  \qquad  {\rm and} 
\qquad  [\![A_0]\!]_2\geq [\![A'']\!]_2\geq \delta \vert A''\vert .
\end{displaymath}
 
Adding those two inequalities gives \eqref{3.7bis}. 

Now, \eqref{3.3} and \eqref{3.7bis} allow us to find  a 
quasi-independent set $B_1 \subset A_0$ such that 
\begin{displaymath}
\vert B_1\vert \geq K^{-1} \frac{[\![A_0]\!]_2^2}{\vert A_0\vert} 
\geq K^{-1} \frac{\delta}{2} \, [\![A_0]\!]_2.
\end{displaymath}

We first notice that we can assume that $\vert B_1\vert\leq K^{-1}\delta [\![A_0]\!]_2$, provided that we reduce $B_1$. Indeed, this can be done as far as we are sure that $K^{-1}\delta [\![A_0]\!]_2\ge1$. This latter fact can be proved in the following way: if we had $1\ge K^{-1}\delta [\![A_0]\!]_2$ then $K\delta^{-1}\ge[\![A_0]\!]_2\ge|A_0|^{1/2}$. Since $[\![A_0]\!]_2\ge c|A_0|^{1/2}(\log|A_0|)^{1/2}$, for some absolute constant $c>0$, and $\delta\ge\displaystyle\frac{[\![A_0]\!]_2}{\vert A_0\vert}$, we would have $\delta\ge c|A_0|^{-1/2}(\log|A_0|)^{1/2}\ge cK^{-1}\delta(\log|A_0|)^{1/2}$. Hence $|A_0|$ would be less than $e^{(Kc^{-1})^2}$. We conclude that $\delta$ would be greater than an absolute constant $e^{-(Kc^{-1})^2}=\delta_0>0$, which is wrong up to a choice of $\delta$ small enough (for instance less than $\delta_0/2$) at the beginning of the proof.

If $\vert B_1\vert\geq \vert A_0 \vert/ 2$, we stop. 

Otherwise, we proceed as follows: suppose more generally 
that we have found disjoint quasi-independent sets  $B_1, \ldots, B_N\subset A_0$ 
satisfying \eqref{3.4}. If they also verify \eqref{3.5}, we stop. If they do not, we set: 
\begin{displaymath}
A' = A_0 \setminus (B_1\cup \cdots\cup B_N),\qquad  \vert A'\vert \geq \frac{\vert A_0\vert}{2} \,\cdot
\end{displaymath}
As before, we can find $B_{N+1}\subset A'$, quasi-independent, such that:
\begin{displaymath}
\vert B_{N+1}\vert \geq  K^{-1} \frac{[\![A']\!]_2^2}{\vert A'\vert} \geq K^{-1}\delta^2 \vert A'\vert
\geq K^{-1}\delta^2\frac{|A_0|}{2}\geq \frac{ K^{-1}\delta}{2}\,  [\![A_0]\!]_2 ,
\end{displaymath}
and we can also assume that $\vert B_{N+1}\vert\leq  K^{-1}\delta  [\![A_0]\!]_2$. Therefore, after a 
finite number of steps, we will have performed \eqref{3.4} and \eqref{3.5}. And, as we already said, this ends 
the proof of Theorem~\ref{theo p-stable stat = Sidon}.\hfill $\square$

\medskip
We can extend Theorem~\ref{extractgauss} from the gaussian to the general $p$-stable framework.

\begin{theorem}\label{extractstable}
There exists a numerical constant $K=K_p$ (depending only on $p$) such that the following inequality holds :
\begin{equation} 
K^{-1} q(A)\leq \Vert i_{A,p}\Vert^{p'}\leq K q(A).
\end{equation} 
\noi In particular, we have :
\begin{equation}
q(A) \geq K^{-1} \bigg(\frac{[\![A]\!]_p}{\;\vert A\vert^{1/p}}\bigg)^{p'}\,\cdot
\end{equation}
\end{theorem}
\smallskip

\noi {\bf Remark.} The lower bound given by this theorem has to be compared to the one of a different kind (involving Orlicz funcions) given by \cite{LeLiQuRo}, Theorem 3.2. Actually, both inequalities are useful as this will be the case in the proof of Theorem~\ref{orlicz}.
\smallskip

\noi{\bf Proof.} The lower bound is easy to obtain: let $B\subset A$ be a quasi-independent set such that $|B|=q(A)$. We have $[\![B]\!]_p\ge c|B|$, for some $c>0$ (depending only on $p$), since $B$ is a Sidon set with a universal constant.

Then, testing the norm of $i_{A,p}$ on the function $f=\displaystyle\sum_{\gamma\in B}\gamma$, we have 
\begin{displaymath}
\Vert i_{A,p}\Vert\ge\frac{[\![B]\!]_p}{\;|B|^{1/p}}\ge c|B|^{\frac{1}{p'}}=c(q(A))^{\frac{1}{p'}}
\end{displaymath}
which was the claim.

For the upper bound, take any $f\in F_p$. There exists a polynomial $P$ such that $\hat P\ge1/4e$ on $A$, $\|P\|_1=1$ and $\log_2\|P\|_\infty\le5eq(A)$ (see  \cite{Ro}, Lema 1.2 or \cite{LiQ}, p.~513). Then $[\![f]\!]_p\le 4e[\![f\ast P]\!]_p$ by the contraction principle.

Thanks to Theorem~\ref{theo Marcus-Pisier}, we have 
\begin{displaymath}
[\![f\ast P]\!]_p\le C\|f\|_{F_p}.\|P\|_{L^{\varphi_{p'}}},
\end{displaymath}
for some $C>0$ (depending on $p$ only).

Finally, 
\begin{displaymath}
\|P\|_{L^{\varphi_{p'}}}\approx \int_G|P|\big(1+\log(1+|P|)\big)^{1/p'}dx\le \|P\|_1.\big(1+\log(1+\|P\|_\infty)\big)^{1/p'}
\end{displaymath}

Hence, we have some $k>0$ (depending on $p$ only) such that 
\begin{displaymath}
\|P\|_{L^{\varphi_{p'}}}\le k\big(q(A)\big)^{1/p'}
\end{displaymath}

The conclusion follows: $[\![f]\!]_p\le K\|f\|_{F_p}.\big(q(A)\big)^{1/p'}$.   \hfill $\square$
\bigskip\goodbreak

\noi {\bf Remark.} Let us emphasize that Proposition \ref{prop Sidon} can be proved very quickly using the previous theorem in the particular case of $p$-stable variables. Indeed, testing the hypothesis of this proposition with $f=\displaystyle\sum_{\gamma\in A}\gamma$, where $A$ is any finite subset of $\Lambda$, we have: $[\![f]\!]_p\ge C^{-1}|A|$ so that $q(A)\ge c|A|^{p'(1-1/p))}=c|A|$ for some constant $c>0$ (depending only on $p$). 

\section{$p$-stable $q$-Rider sets}

Our aim in this Section is to introduce an apparently new notion of thin set, that of  \emph{$p$-stable 
$q$-Rider set}, and to compare it with the previously known notion of \emph{$q$-Rider set}, which might 
be called \emph{$2$-stable $q$-Rider set} to recall  that it is defined with Gaussian (equivalently Rademacher) 
variables. We will always assume that $p,q$ are given in such a way that  
\begin{displaymath}
1\leq q < p\leq 2,  
\end{displaymath}
and we will say that $\Lambda\subset\Gamma$ is a \emph{$p$-stable $q$-Rider set} if the following a priori inequality holds:
\begin{equation}\label{4.1} 
\qquad \Vert f \Vert_{F_q}\leq C\, [\![ f ]\!]_p, \qquad  \forall f\in {\cal P}_\Lambda.
\end{equation}

The reader should compare with \eqref{1.7} to see what is new here. It should be emphasized that the 
behaviour of the $[\![\  ]\!]_p$ -norm is very different from that of the $[\![\  ]\!]_2$-one: for example, if 
$G =\T$,  $e_{n} (t) = \e^{int}$ and $f (t) =\sum_{n=1}^N a_n e_{n} (t)$, Marcus and Pisier  
(\cite{MaPi}, Remark~5.7), extending a result of Salem and Zygmund  for $p = p' = 2$ (see 
\cite{LiQ}, Chapitre~13, Proposition~III.13), proved that, if $a = (a_1, \ldots, a_N)$ and $N\geq 2$: 
\begin{equation}\label{4.2} 
[\![ f ]\!]_p \geq C_0\, N^{1/p}(\log N)^{1/p'}\,\frac{1}{N} \sum_{n=1}^N |a_n|
\end{equation} 
and that the reverse inequality holds for $a_1= \cdots = a_N = 1$. So, we might  expect that the new 
notion introduced is highly depending on $p$ (and on $q$, of course\,!). This is not quite the case, as will be 
apparent in our results, which will mainly consist in two theorems:  a ``functional-type'' condition, indicating 
that those sets can be defined by other a priori inequalities, and an \emph{equivalence Theorem} showing that 
indeed this notion is nothing but $s$-Riderness for some value of $s$, depending on $p$ and $q$. 
It will be convenient to recall the following definitions and facts:
\begin{enumerate}
\item For $r > 0$, $\psi_{r}$ will denote the Orlicz function $x\mapsto \e^{x^{r}} - 1$, $x \geq 0$. \\
If $A\subset \Gamma$ is finite, we set 
\begin{equation}\label{psi indice A}
\psi_{r}(A) =\Big\Vert\sum_{\gamma\in A} \gamma \,\Big\Vert_{\psi_r}.
\end{equation}
\item The space $F_p$ and its norm have already been defined in the Introduction: see \eqref{1.6}.
\item If $X,Y$ are two Banach spaces continuously contained in $L^1 (G)$, the set ${\cal M} (X,Y)$ of 
multipliers of $X$ to $Y$ is the set of families $m = (m_\gamma)_{\gamma\in\Gamma}$ of complex 
numbers  such that, whenever  $f =\sum_\gamma a_\gamma\gamma\in X$, then 
\begin{displaymath}
\qquad \qquad
g = \sum_\gamma m_\gamma  a_\gamma\gamma \in Y, \qquad 
{\rm with} \quad \Vert g\Vert_Y\leq C\Vert f\Vert_X.
\end{displaymath}
The best constant $C$ being called the multiplier-norm of $m$ and being denoted by 
$\Vert m\Vert_{{\cal M}(X,Y)}$.
 \item We denote by ${\cal C}^{p-as}$ the completion of ${\cal P}$ with respect to the $[\![ \,\  ]\!]_p$-norm; 
it is the so-called Banach space of almost surely continuous $p$-stable random Fourier series (\cite{MaPi}). 
Then, the dual space of  ${\cal C}^{p-as}$ is isomorphic to ${\cal M}(F_p, L^{\psi_{p'}})$, the set of 
multipliers from $F_p$ to $L^{\psi_{p'}}$. This  result  (well-known for $p=2$: see \cite{MarPis})  follows in 
a standard way from the delicate Theorem~\ref{theo Marcus-Pisier}, as it is already the case for $p=2$ and we 
will not detail that formal  extension. 
\item Once and for all, we set 
\begin{align}
\varepsilon &= \frac{p - q}{q(p-1)} \label{epsilon}=1-\frac{p'}{q'}\cdot 
\\
\frac{1}{\alpha} & = \frac{1}{p}+ \frac{1}{q'} \label{alpha}
\\
\beta & = \frac{\eps}{p'} + \frac{1}{p} - \frac{1}{2}= \frac{1}{q} - \frac{1}{2}\cdot\label{beta}
\end{align}
\end{enumerate}

\medskip

We first  prove the  following simple  proposition (the case $q=2$ being already known \cite{LiQuRo}), 
which will actually follow from Theorem~\ref{theo equivalences} below, but which will motivate this theorem.

\begin{proposition} [{\bf Mesh condition}] \label{mesh condition}
Let $\Lambda\in\Z$ be a $p$-stable $q$-Rider set. Then, there exists a constant $K$ such that, for each integer 
$N\geq 2$, one has :
\begin{equation}\label{4.4}
\vert\Lambda\cap\{1, \ldots, N\}\vert\leq K (\log N)^{\frac{(p - 1)q}{p-q}} 
= K (\log N)^{\frac{1}{\varepsilon}}.
\end{equation}
\end{proposition}

\noi {\bf Proof.} Set $B=\Lambda\cap\{1, \ldots, N\} =\{\lambda_1 < \cdots <\lambda_n\}$, and 
$f =\sum_{j=1}^n e_{\lambda_j}$. The assumption and Lemma~\ref{majoriz. p-stable norm} give us, for 
some constant $C$, that 
$n^{\frac{1}{q}} \leq C [\![ f ]\!]_p\leq C C_0 n^{\frac{1}{p}} (\log N)^{\frac{1}{p'}}$. Grouping terms gives 
$n\leq K (\log N)^{\frac{(p - 1)q}{p - q}}$,  where we  set $K = (C C_0)^{\frac{pq}{p - q}}$, proving our 
claim.\hfill $\square$
\bigskip

We will now prove the main result of this section: 
\medskip

\begin{theorem}\label{stableqrider}
Let $\Lambda\subset\Gamma$.
\vskip 5pt
\noindent $\Lambda$ is a $p$-stable $q$-Rider set if and only if $\Lambda$ is an $s$-Rider set, with $s =\frac{2q'}{2q'-p'}\,\cdot$
\end{theorem}

\medskip

Actually we are going to prove the following more precise theorem:
\medskip
\begin{theorem}\label{theo equivalences}
Let $\Lambda\subset\Gamma$. Then, the following conditions are equivalent:\par
\vskip 2pt
$(1)$ $\Lambda$ is a $p$-stable $q$-Rider set;\par\vskip 2pt
$(2)$ $\ell_{q'} (\Lambda) \hookrightarrow  {\cal M}(F_p, L^{\psi_{p'}})$;\par\vskip 2pt
$(3)$ $\ell_{\alpha} (\Lambda)\hookrightarrow  L^{\psi_{p'}}$;\par\vskip 2pt
$(4)$ $\psi_{p'} (A) \leq C\,|A|^{1/\alpha}$ for some constant $C >0$, and every finite 
subset $A$ of $\Lambda$;\par\vskip 2pt
$(5)$ $q(A) \geq c\,\vert A\vert^\varepsilon$, for some constant $c > 0$, and for all finite subsets 
$\{0\} \not= A\subset \Lambda$;\par\vskip 2pt
$(6)$ $\Lambda$ is an $s$-Rider set, with $s =\frac{2q'}{2q'-p'}= \frac{2q (p - 1)}{p - 2q + pq}\,\cdot$
\end{theorem}

Point out the simple relation: $2q'=s'p'$. Moreover, to precise the behavior for the ``degenerate'' cases: when $q=1$, we have $s=1$ as well (remember Proposition \ref{prop Sidon}). On the other hand, when $s=1$, $q=1$ and any $p$ fits.

\medskip
Let us make some comments. This result is known for $p=2$ and has been proved by the fourth-named 
author (\cite{Ro2}, \cite{Ro}, Teorema~III.2.3). The symbol $\hookrightarrow $ means that the left-hand space 
is mapped to the right-hand one by means of the Fourier transform or of its inverse. Recall that $q(A)$ denotes 
the largest possible cardinality of a quasi-independent subset of $A$, and that the definition of 
quasi-independent sets is given at the beginning of the proof of  Theorem~\ref{theo p-stable stat = Sidon}. \par
\medskip

\noi {\bf Remark 1.} We know (see \cite{LiQuRo}) that the \emph{mesh condition} for $s$-Rider reads  as 
$\vert\Lambda\cap\{1,\ldots, N\}\vert\leq K (\log N)^{\frac{s}{2 - s}}$. But, 
$\frac{s}{2 - s} = \frac{1}{\varepsilon}$, and we fall 
again on \eqref{4.4} of Proposition~\ref{mesh condition}.
\par\bigskip\goodbreak

\noi {\bf Remark 2.} The preceding theorem can be read in two ways. First, any $p$-stable $q$-Rider set is actually an $s$-Rider set for the right value of $s$. But on the other hand, if one fixes some $s\in(1,2)$ and a set $\Lambda$ which is an $s$-Rider set, one can choose either $p$ or $q$ in order to realize $\Lambda$ as a $p$-stable $q$-Rider set. Nevertheless, one has to be careful. Let us precise this:

\hskip1cm a. If one fixes $p\in(1,2]$, then one can choose $q$ such that $q'=s'p'/2$ and $\Lambda$ is a $p$-stable $q$-Rider set. Point out that $q<p\le2$.

\hskip1cm b. If one fixes $q\in(1,2)$, then one can choose $p\in(1,2]$ such that $p'=\frac{2q'}{s'}$ if and only if we have $s\ge q$. In that case,  $\Lambda$ is a $p$-stable $q$-Rider set.

\par\bigskip\goodbreak

\noi {\bf Proof of Theorem~\ref{theo equivalences}.} \par
\noi $(1) \Leftrightarrow (2)$. The Fourier transform maps $X_{\Lambda}$ to $\ell_{q}(\Lambda)$ if and only if its 
transpose maps $\ell_{q'}(\Lambda)$ to the dual of $X_{\Lambda}$. The result now easily follows from 
the previous description of the dual of $X$.\par

\smallskip

\noi $(2) \Rightarrow (3)$. Let $f =\displaystyle\sum_{\gamma \in\Lambda} c_\gamma\gamma\in {\cal P}_\Lambda$. 
We can write (in short: see \eqref{alpha}, we have $\ell_\alpha = \ell_p \, . \, \ell_{q'}$):
\begin{displaymath}
c_\gamma=a_\gamma b_\gamma,
\end{displaymath}
where
\begin{displaymath}
\quad \vert a_\gamma\vert =\vert c_\gamma\vert^{1-\theta}, \qquad  
\vert b_\gamma\vert =\vert c_\gamma\vert^{\theta},
\end{displaymath}
with $\theta = \frac{\alpha}{p} = 1 - \frac{\alpha}{q'}$ and
\begin{displaymath}
\quad \Vert a\Vert _{q'} = \Vert c \Vert_\alpha^{1 -\theta}, \qquad  \Vert b\Vert _{p} 
= \Vert c\Vert_\alpha^{\theta}.
\end{displaymath}

If we set $P=\sum_\gamma b_\gamma\gamma$, we have from $(2)$ that 
$\Vert f \Vert_{\psi_{p'}} \leq C\,\Vert a\Vert_{q'} \Vert P \Vert _{F_p}$, where $C$ is some constant. 
Equivalently: 
$\Vert f \Vert_{\psi_{p'}} \leq C\, \Vert c\Vert_\alpha^{1-\theta +\theta} = C \, \Vert c\Vert_\alpha$, 
which was our claim.\par

\smallskip

\noi $(3) \Rightarrow (2)$. This is obvious since $\frac{1}{\alpha}= \frac{1}{p}+ \frac{1}{q'}\cdot$
\smallskip

\noi $(3) \Rightarrow (4)$. Indeed, 
$\psi_{p'}(A) \leq C\, \Vert\widehat{\ind_A} \Vert_\alpha = C\,\vert A\vert^{\frac{1}{\alpha}}$. 
\smallskip

\noi $(4) \Rightarrow (5)$. We use a result that we proved in (\cite{LeLiQuRo}, Proposition~3.2), namely that, 
for any finite set  $A\subset \Gamma$, $A \not= \{0\}$, and any $r > 0$, we have: 
\begin{equation}\label{4.5} 
q(A)\geq C_r \Big(\frac{\vert A\vert}{\psi_r (A)} \Big)^r
\end{equation}
We use \eqref{4.5} with $r = p'$. Our assumption implies that 
$q(A) \geq C \vert A\vert^{(1- \frac{1}{\alpha}) p'}$, and
\begin{displaymath}
\Big(1 - \frac{1}{\alpha} \Big) p' = \Big(1 - \frac{1}{q'} - \frac{1}{p} \Big) p' 
= \Big(\frac{1}{p'} - \frac{1}{q'} \Big) p'=1- \frac{p'}{q'}=\varepsilon.
\end{displaymath}
\smallskip

\noi $(5) \Rightarrow (6)$. By \cite{Ro2} or \cite{Ro}, Teorema~III.2.3, we know that this condition is equivalent to the fact that $\Lambda$ is an $s$-Rider set where $\varepsilon= \frac{2}{s} - 1$. As $\varepsilon=1-\frac{p'}{q'}$, we have $s=\frac{2q'}{2q'-p'}$.
\smallskip

\noi $(6) \Rightarrow (3)$. Either using \cite{Ro2}, \cite{Ro}, Teorema~III.2.3 or applying the equivalence of $(1)$ and $(3)$ when $p=2$ and $s$ instead of $q$, we know that $\ell_{\tilde\alpha} (\Lambda)\hookrightarrow  L^{\psi_{2}}$ where $\frac{1}{\tilde\alpha}=\frac{1}{2}+\frac{1}{s'}$. On the other hand, we obviously have $\ell_{1} (\Lambda)\hookrightarrow  L^\infty$. As $p>q$, we have $\tilde\alpha>\alpha$, so a standard interpolation argument implies that $\ell_{\alpha} (\Lambda)\hookrightarrow  L^{\psi_{r}}$ with 

\begin{displaymath}
\frac{1}{\alpha} =\frac{\theta}{\tilde\alpha}+\frac{1-\theta}{1} \hskip1cm\hbox{and}\hskip1cm \frac{1}{r}=\frac{\theta}{2}+\frac{1-\theta}{\infty}
\end{displaymath}
We obtain $\theta=\displaystyle\Big(\frac{1}{q'} - \frac{1}{p'}\Big)\Big(\frac{1}{s'} - \frac{1}{2} \Big)^{-1}=\frac{2}{p'}$, since $s'p'=2q'$. Hence  $r=\frac{2}{\theta}=p'$.\hfill $\square$
\bigskip

The following theorem has two aims. First the proof of the preceding theorem is not self-contained: to prove that $(5) \Rightarrow (6)$, we used the fact that the theorem was already known for $p=2$. The following theorem provides a proof of this result as well.

On the other hand, this will actually provide a stronger result, which cannot a priori  be obtained just using the case $p=2$ (i.e. assuming the results of \cite{Ro2} or \cite{Ro}). Nevertheless, the proof proceeds as in \cite{Ro}, using a difficult lemma of Bourgain on 
quasi-independent sets. 

Though there is essentially no change with regard to \cite{Ro}, we will give the details, for the convenience of the reader. The links beetween the values of the parameters are the same as before.\par

\begin{theorem}\label{Lorentz}
Let $\Lambda\subset \Gamma$. The following conditions are equivalent:\par
\vskip 2pt

$(i)$ $q(A) \geq c\,\vert A\vert^\varepsilon$, for some constant $c > 0$, and for all finite subsets 
$\{0\} \not= A\subset \Lambda$;\par\vskip 2pt

$(ii)$  ${\cal C}^{p-as}_{\Lambda} :=  X_{\Lambda} \hookrightarrow \ell_{q,1}(\Lambda)$;\par\vskip 2pt

$(iii)$ $\Lambda$ is a $p$-stable $q$-Rider set;\par\vskip 2pt

$(iv)$ ${\cal C}^{p-as}_{\Lambda} :=  X_{\Lambda} \hookrightarrow \ell_{q,\infty}(\Lambda)$.

\end{theorem}

Recall that $\ell_{q,1} = \ell_{q,1} (\Lambda)$ (resp. $\ell_{q,\infty} = \ell_{q,\infty} (\Lambda)$) is the Lorentz space of families 
$a = (a_\gamma)_{\gamma\in\Lambda}$ tending to $0$ whose decreasing 
rearrangement $(a_n^*)_n$ satisfies $\displaystyle\sum_{n\geq 1} \frac{a_n^*}{n^{1/q'}}  <\infty$ (resp. $\displaystyle\sup_{n\geq 1} n^{1/q}a_n^*  <\infty$). 

One has 
$\ell_{q,1}(\Lambda) \hookrightarrow \ell_q (\Lambda) \hookrightarrow\ell_{q,\infty}(\Lambda)$.
\par\medskip

\noi {\bf Proof.} The implications $(ii) \Rightarrow (iii)\Rightarrow (iv)$  are obvious.

$(iv) \Rightarrow(i)$ is easy with help of Theorem \ref{extractstable}: for any finite subset $A$ of $\Lambda$, we have 

\begin{displaymath}
q(A)\geq K^{-1} \bigg(\frac{[\![A]\!]_p}{\vert A\vert^{1/p}}\bigg)^{p'}\ge K_1\big(|A|^{1/q-1/p}\big)^{p'}
\end{displaymath}
which gives $(iv)$, since $p'\big(\frac{1}{q}-\frac{1}{p}\big)=p'\big(\frac{1}{p'}-\frac{1}{q'}\big)=\eps$.

It remains to prove the difficult part ($(i) \Rightarrow (ii)$) of this theorem. We first need the following lemma.

\begin{lemma}[\cite{Ro}, Lema III.2.6]\label{lema thesis}
Let $\Lambda$ be a set satisfying the condition $(5)$ of Theorem~\ref{theo equivalences}. For any finite 
subset $A$ of $\Lambda$ such that $A \not = \{0\}$ and $c\,|A|^\eps \geq 2$, there exist $N$ 
pairwise disjoint quasi-independent sets $B_1, \ldots, B_N \subset A$ such that:\par
$a)$ $\frac{2}{c}\,|A|^{1 - \eps} \geq N \geq \frac{1}{2c} \,|A|^{1 - \eps}$;\par
\vskip 2pt
$b)$ $c\,|A|^\eps \geq |B_j| \geq  \frac{c}{2}\, |A|^\eps$ for every $j=1,\ldots, N$.
\end{lemma}

\noi{\bf Proof of Lemma~\ref{lema thesis}.} By our assumption, there is a quasi-independent subset $B$ of $A$ 
such that $|B| \geq c\,|A|^\eps$. Since $\frac{c}{2}\,|A|^\eps \geq 1$, we can find a (quasi-independent) 
subset $B_1$ of $B$ such that $\frac{c}{2}\,|A|^\eps \leq |B_1| \leq c\,|A|^\eps$.\par
Assume now that pairwise disjoint quasi-independent subsets $B_1, \ldots, B_n$ of $A$ has been constructed 
such that $c\,|A|^\eps \geq |B_j|\geq \frac{c}{2}\, |A|^\eps$ for $1 \leq j \leq n$. There are two 
possibilities:\par
$(i)$ If $|A|/2 > \big| \bigcup_{j=1}^n B_j \big|$, we choose in $A \setminus \bigcup_{j=1}^n B_j$ a 
quasi-independent subset $B_{n+1}$ whose cardinal is
\begin{displaymath}
|B_{n+1} | \geq c\, \Big| A \setminus \bigcup_{j=1}^n B_j \Big|^\eps 
\geq c\, \Big( \frac{|A|}{2} \Big)^\eps \geq \frac{c}{2}\, |A|^\eps,
\end{displaymath}
and which we can also choose such that $|B_{n+1}|\leq c\,|A|^\eps$.\par
$(ii)$ If  $|A|/2 \leq \big| \bigcup_{j=1}^n B_j \big|$, we stop the process at $N=n$. Indeed, we then have:
\begin{displaymath}
\frac{|A|}{2} \leq \Big| \bigcup_{j=1}^N B_j \Big| \leq N\,c\, |A|^\eps
\end{displaymath}
and hence $N \geq \frac{1}{2c}\, |A|^{1 - \eps}$. On the other hand, since $B_1, \ldots, B_N$ are disjoint, 
one has
\begin{displaymath}
|A| \geq \Big| \bigcup_{j=1}^N B_j \Big| = \sum_{j=1}^N |B_j| \geq N\,\frac{c}{2}\,|A|^\eps ;
\end{displaymath}
hence $N \leq \frac{2}{c}\,|A|^{1 - \eps}$.\par
\smallskip

That ends the proof of Lemma~\ref{lema thesis}. \hfill $\square$
\par\medskip

Now, we have to show that there exists a constant $C > 0$ such that
\begin{displaymath}
\| \widehat f \|_{q,1} \leq C\, [\![ f ]\!]_p
\end{displaymath}
for every trigonometric polynomial $f\in{\cal P}_\Lambda$.\par
By homogeneity, we may assume that $\|\hat f\|_\infty =1$, and we consider the level sets
\begin{displaymath}
\qquad A_j = \{ \gamma \in \Lambda\,;\ 2^{- j + 1} \geq | \widehat f (\gamma) | > 2^{-j} \}\,, 
\qquad j=1,2, \ldots
\end{displaymath}

Recall now the following result of Bourgain (\cite{Bourgain1}, Lemma 2; see also \cite{Bourgain2}, and 
\cite{LiQ}, Chapitre~12, Lemme~I.10).

\begin{lemma}[{\bf Bourgain}] \label{lemme Bourgain}
There exists a numerical constant $R > 10$ such that for every family $B_1, \ldots, B_L$ of pairwise disjoint 
finite quasi-independent sets such that
\begin{displaymath}
\qquad \qquad \frac{|B_{l+1} |}{|B_l|} \geq R \,, \qquad \text{for } l = 1, \ldots, L - 1,
\end{displaymath}
one can find, for each $l = 1, \ldots, L$, a subset $C_l \subset B_l$ such that:
\begin{displaymath}
\quad\qquad |C_l| \geq \frac{1}{10}\,|B_l|, \qquad \text{for all } l = 1, \ldots, L,
\end{displaymath}
and the union $\bigcup_{l=1}^L C_l $ is quasi-independent.
\end{lemma}

Setting $R_1 = \max \{R^{1/2\beta}, (2R)^\eps\}$, where $R$ is the above constant and $\beta$ is given by (\ref{beta}), we define $j_1 = 1$ and 
\begin{equation}\label{def j_l}
j_{l+1} = \min \{ j > j_l\,;\ |A_j| > R_1\, |A_{j_l}| \},
\end{equation}
whenever this last set is nonempty; we stop and take $L=l$ when it is empty (this eventually happens since 
$f$ is a trigonometric polynomial). Set
\begin{displaymath}
N_l = |A_{j_l}|.
\end{displaymath}

We have the following upper estimate:
\begin{displaymath}
\| \widehat f\|_{q,1} \leq \sum_{j \geq 1} 2^{-j + 1} \| \ind_{A_j} \|_{q,1} 
\leq \sum_{l=1}^L \sum_{j_l \leq j < j_{l+1}} 2^{-j +1} q\, |A_j|^{1/q}\,,
\end{displaymath}
since
\begin{displaymath}
\| \ind_{A} \|_{q,1}  = \sum_{n=1}^{|A|} \frac{1}{n^{1/q'}} 
\leq \int_0^{|A|} \frac{dx}{x^{1/q'}} = q\, |A|^{1/q}.
\end{displaymath}
Now, we have $|A_j| \leq R_1 |A_{j_l}|$ for $j_l \leq j < j_{l+1}$; hence we get
\begin{align}
\| \widehat f\|_{q,1} 
& \leq q\sum_{l=1}^L R_1^{1/q} 2^{-j_l} |A_{j_l}|^{1/q} \sum_{j_l \leq j < j_{l+1}} 2^{-j +1 + j_l} \notag \\
& \leq 4q R_1^{1/q} \sum_{l=1}^L 2^{-j_l} N_l^{1/q}. \label{upper estimate}
\end{align}

We are now going to use Lemma~\ref{lema thesis}. For this purpose, let us denote by $T$ the first index 
$l$ such that $c\,N_l^{2\beta} \geq 2$. If no such an index exists, we will set $T = L + 1$.\par
We will split the sum \eqref{upper estimate} into two parts.\par\smallskip

First, the Cauchy-Schwarz inequality gives
\begin{align*}
\sum_{l=1}^{T-1} 2^{-j_l} N_l^{1/q} 
& \leq \Big( \sum_{l=1}^{T-1} 2^{-2j_l} N_l \Big)^{1/2} 
\Big( \sum_{l=1}^{T-1} N_l^{2\beta} \Big)^{1/2} \\
& \leq \Big(\sum_{\gamma \in \Lambda} |\widehat f (\gamma) |^2 \Big)^{1/2} 
\bigg[ \sum_{l=1}^{T-1} \Big(\frac{N_{T-1}}{R_1^{T - l -1}} \Big)^{2\beta} \bigg]^{1/2} \\
& = \| f \|_2 N_{T-1}^\beta \Big(\sum_{k \geq 0} R_1^{-2k \beta} \Big)^{1/2} .
\end{align*}
Since $N_{T-1}^{2\beta} < 2/c$ and $R_1^{2\beta} \geq R > 10 > 4$, we get 
\begin{displaymath}
\sum_{l=1}^{T-1} 2^{-j_l} N_l^{1/q} \leq \| f \|_2 \Big(\frac{2}{c}\Big)^{1/2} \Big(\frac{4}{3} \Big)^{1/2}
\end{displaymath}
and
\begin{equation}\label{piece 1}
\sum_{l=1}^{T-1} 2^{-j_l} N_l^{1/q} \leq \frac{2}{\sqrt c}\, [\![ f ]\!]_p .
\end{equation}
\medskip

For $l \geq T$, we apply Lemma~\ref{lema thesis} to the set $A_{j_l}$: we get $M_l$ pairwise disjoint 
quasi-independent subsets $B_{l,1}, \ldots, B_{l, M_l} \subset A_{j_l}$ such that
\begin{equation}\label{size sets B}
\qquad \quad \frac{c}{2} N_l^{\eps} \leq |B_{l,m} | \leq c N_l^\eps\,, \qquad m= 1, \ldots, M_l,
\end{equation}
with
\begin{equation}\label{size number M}
\qquad \quad \frac{1}{2c} \, N_l^{1 -\eps} \leq M_l \leq \frac{2}{c}\,N_l^{1 - \eps}\,, \qquad 
l= T, \ldots, L.
\end{equation}
Since $N_{l+1} \geq R_1 N_l$, we get, for $T \leq l < L$:
\begin{equation}\label{ratios M}
\frac{M_L}{M_l} \geq \frac{R_1^{1 - \eps}}{4} \geq \frac{1}{4}
\end{equation}
and
\begin{equation}\label{ratios B}
\frac{|B_{l+1,m} |}{|B_{l,m'}|} \geq \frac{R_1^\eps}{2} \geq R. 
\end{equation}

Using \eqref{ratios M}, we can find, for each $l = T, \ldots, L - 1$, a map 
\begin{displaymath}
\phi_l \colon \{1, \ldots, M_L\} \to \{1, \ldots, M_l\} 
\end{displaymath}
such that 
\begin{displaymath}
\qquad \quad | \phi_l^{-1} (m) | \leq 4\,\frac{M_L}{M_l}_,,\qquad m= 1, \ldots, M_l.
\end{displaymath}
Applying, for each $m=1, \ldots, M_L$,  Lemma~\ref{lemme Bourgain} to the 
sequence
\begin{displaymath}
B_{T,\phi_T (m)}, \ldots, B_{L-1, \phi_{L-1} (m)}, B_{L,m} \,,
\end{displaymath}
thanks to \eqref{ratios B}, we get, for each $l = T, \ldots, L$ and each $m = 1, \ldots, M_L$, a 
quasi-independent set $C_{l,m}$ such that
\begin{equation}\label{size sets C}
\hskip 30 pt | C_{l,m}| > \frac{1}{10} \, |B_{l, \phi_l (m)} | \geq \frac{c}{20}\, N_l^\eps \,, 
\quad l=T, \ldots, L,
\end{equation}
and such that
\begin{equation}\label{qi}
\bigcup_{l=T}^L C_{l,m} \text{ is quasi-independent for } m= 1, \ldots , M_L.
\end{equation}

We now introduce, for $m=1, \ldots, M_L$, the trigonometric polynomial
\begin{displaymath}
g_m = \sum_{l=T}^L \Big( \frac{M_l}{4 M_L}\Big)^{1/p} 2^{-j_l} \sum_{\gamma \in C_{l,m}} \gamma.
\end{displaymath}

We claim that, for every $\gamma\in \Gamma$, one has
\begin{equation}\label{condition lower p-estimate}
| \widehat f (\gamma) |^p \geq \sum_{m=1}^{M_L} | \widehat g_m (\gamma) |^p.
\end{equation}

Indeed, it suffices to check that for $\gamma$ in the spectra of the $g_m$'s, and if $\gamma$ is in some 
set $B_{l,m_0}$, it cannot be in another one, because these sets are pairwise disjoint; hence
\begin{align*}
\sum_{m=1}^{M_L} | \widehat g_m (\gamma) |^p 
& = 2^{-p j_l} \frac{M_l}{4 M_L}\,\big| \{m\,;\ \gamma \in C_{l,m} \} \big| 
\leq 2^{-p j_l} \frac{M_l}{4 M_L}\, | \phi_l^{-1} (m_0) | \\
& \leq 2^{-p j_l}  \leq | \widehat f (\gamma) |^p.
\end{align*}

It follows from the lower $p$-estimate of the norm $[\![ \,\, ]\!]_p$ (Theorem~\ref{lower p-estimate}) that
\begin{equation}\label{lower p-estimate of f}
a\,[\![ f ]\!]_p \geq \Big(\sum_{m=1}^{M_L} [\![ g_m ]\!]_p\Big)^{1/p}.
\end{equation}

But, by \eqref{qi}, the spectrum of $g_m$ is quasi-independent, and every quasi-independent set is a Sidon set, 
with constant  $\leq 5$ (see the beginning of the proof of  Theorem~\ref{theo p-stable stat = Sidon}), so 
$\| \widehat g_m \|_1 \leq 5\,\|g_m\|_\infty$, and $\| \widehat g_m \|_1 \leq 5\,[\![ g_m]\!]_p$, 
by Proposition~\ref{converse Sidon}. It follows, from \eqref{size number M} and \eqref{size sets C}, that:
\begin{align*}
5\,[\![ g_m]\!]_p 
& \geq \sum_{l=T}^L | C_{l,m}| \Big( \frac{M_l}{4 M_L} \Big)^{1/p} 2^{-j_l} 
\geq \frac{1}{4^{1/p} M_L^{1/p}} \sum_{l=T}^L \frac{c}{20} \, N_l^\eps 
\Big( \frac{N_l^{1 - \eps}}{2c}\Big)^{1/p} 2^{-j_l} \\
& = \frac{c^{1/p'}}{20(8M_L)^{1/p}} \sum_{l=T}^L N_l^{\eps + \frac{1 - \eps}{p}} 2^{-j_l} .
\end{align*}
Since 
\begin{displaymath}
\eps + \frac{1 - \eps}{p} = \big( 1 - \frac{1}{p} \Big) \eps + \frac{1}{p} = \frac{\eps}{p'} + \frac{1}{p} 
= \Big( \frac{1}{q} - \frac{1}{p} \Big) + \frac{1}{p} = \frac{1}{q} \,,
\end{displaymath}
we get
\begin{equation}\label{estim g_m}
[\![ g_m]\!]_p \geq 
\frac{c^{1/p'}}{100(8M_L)^{1/p}} \sum_{l=T}^L N_l^{1/q} 2^{-j_l} .
\end{equation}

Therefore \eqref{lower p-estimate of f} gives
\begin{align}
[\![ f ]\!]_p 
& \geq \frac{c^{1/p'}}{100 a (8M_L)^{1/p}} \,\bigg[\sum_{m=1}^{M_L} 
\Big(\sum_{l=T}^L N_l^{1/q} 2^{-j_l} \Big)^p \bigg]^{1/p} \notag\\
& = \frac{c^{1/p'}}{100 a . 8^{1/p}} \sum_{l=T}^L N_l^{1/q} 2^{-j_l} . \label{second piece}
\end{align}
\medskip

Putting \eqref{second piece} together with \eqref{piece 1}, we get:
\begin{displaymath}
\sum_{l=1}^L N_l^{1/q} 2^{-j_l} \leq 
\Big( \frac{2}{\sqrt c} + \frac{100 a . 8^{1/p}}{c^{1/p'}} \Big) \, [\![ f ]\!]_p.
\end{displaymath}

It remains to use \eqref{upper estimate} to obtain:
\begin{displaymath}
\| \widehat f \|_{q,1} \leq 
4q R_1^{1/q} \,\Big( \frac{2}{\sqrt c} + \frac{100 a . 8^{1/p}}{c^{1/p'}} \Big) \, [\![ f ]\!]_p,
\end{displaymath}
and achieve the proof of Theorem~\ref{Lorentz}. 
\hfill $\square$


\section{The link with Orlicz spaces}

We are going to characterize $s$-Rider sets in terms of continuous mapping to Orlicz spaces (remember the beginning of Section 2). 
\medskip
\begin{theorem}\label{orlicz}
Let $\Lambda\subset \Gamma$, $s\in(1,2)$ and $r$ be greater than both $2$ and $\rho=\frac{2-s}{s-1}$. Let $\tilde p=\frac{2r}{2r-\rho}\cdot$
\vskip 2pt
The following conditions are equivalent:\par
\vskip 2pt

$(i)$ $\Lambda$ is an $s$-Rider set;\par\vskip 2pt

$(ii)$  ${\cal C}^{\tilde p-as}_{\Lambda}\hookrightarrow L^{\psi_r}$;\par\vskip 2pt

$(iii)$  For every finite subset $A$ of $\Lambda$, we have $\psi_r(A)\le C [\![A ]\!]_{\tilde p}$, where $C$ does not depend on $A$.
\end{theorem}
\medskip

\noindent{\bf Proof.} $(i)\Rightarrow(ii)$. We already know that we can realize $\Lambda$ as a $p$-stable $q$-Rider set with $p=r'\le2$. Then the value of $q$ is fixed by the relation $q'=rs'/2$. By $(3)$ of Theorem~\ref{theo equivalences}, we know that $\ell_{\alpha} (\Lambda)\hookrightarrow  L^{\psi_{r}}$ with $\frac{1}{\alpha}= \frac{1}{p}+ \frac{1}{q'}=1-\frac{2-s}{sr}$. Let us point out that $\alpha<2$.

 Now we can use Theorem~\ref{theo equivalences} again to realize $\Lambda$ as a $\tilde p$-stable $\alpha$-Rider set but only (see Remark 2 after that theorem) when $s\ge\alpha$. This condition is fulfilled since it is equivalent to the condition $r\ge\rho$. 

Moreover the value of $\tilde p$ is fixed by the relation $\tilde p{\,'}=\frac{2\alpha'}{s'}=\frac{2r}{2r-\rho}\cdot$ The conclusion follows.

$(ii)\Rightarrow(iii)$. This is obvious.

$(iii)\Rightarrow(i)$. Fix any finite subset $A$ of $\Lambda$. We shall use $\eps=\frac{2}{s}-1\cdot$

If $[\![A ]\!]_{\tilde p}\le|A|^{1-\frac{\eps}{r}}$. Then $\psi_r(A)\le C|A|^{1-\frac{\eps}{r}}$. Hence, by  \cite{LeLiQuRo}, Proposition 3.2, 

\begin{displaymath}
q(A)\ge c_r \bigg(\frac{|A|}{\psi_r(A)}\bigg)^r\ge C'_r|A|^{\eps}.
\end{displaymath}
\smallskip

If not, then $[\![A ]\!]_{\tilde p}\ge|A|^{1-\frac{\eps}{r}}$; but

\begin{displaymath}
q(A) \geq K^{-1} \bigg(\frac{[\![A]\!]_{\tilde p}}{\;\vert A\vert^{1/{\tilde p}}}\bigg)^{{\tilde p}'}\,
\end{displaymath}
by Theorem~\ref{extractstable}, so we obtain that $q(A)\ge c|A|^{(1-\frac{\eps}{r}-\frac{1}{\tilde p}){\tilde p}'}$. 

Now, a quick computation gives ${\tilde p}{\,'}=\frac{2r}{\rho}$ and we conclude that

\begin{displaymath}
q(A)\ge c|A|^{\eps}.
\end{displaymath}

So, in every case, we have $q(A)\ge c|A|^{\eps}$ and this characterizes the fact that $\Lambda$ is an $s$-Rider set.\hfill $\square$
\medskip

\noi {\bf Remark.} The preceding theorem extend Theorem 3.1 of \cite{LeLiQuRo}. More precisely, when $s\le\frac{4}{3}$, we can choose $r=\rho$, so $\tilde p=2$ and we recover the version of \cite{LeLiQuRo}.

When $s\ge\frac{4}{3}$, we can take $r=2$ and this gives $\tilde p=\frac{4(s-1)}{5s-6}\cdot$
\bigskip

The previous result leads naturally to investigate more specifically thin sets involving both random norms and Orlicz spaces. This was done for instance by the authors in \cite{LeLiQuRo}, where (among other things) they studied the notion of $\Lambda^{as} (q)$-sets. This latter notion is actually weaker than of the notion of $s$-Rider, or equivalently (through the previous theorem) the notion of what we could call $\Lambda^{p-as} (\psi_r)$-set. It is not known whether the notion of $\Lambda^{as} (q)$-sets is actually different of the usual notion of $\Lambda (q)$-sets. 

In the following, we add several results on $\Lambda^{as} (q)$-sets. Let us first precise the definition.

%
\begin{definition}
A subset $\Lambda$ of $\Gamma$ is said to be a $\Lambda^{p-as} (q)$-set, $1 < p \leq 2$, $q > 2$, 
if there is a constant $C > 0$ such that, for every trigonometric polynomial $f \in {\cal P}_\Lambda$, one has:
\begin{displaymath}
\| f \|_q \leq C\, [\![ f ]\!]_p.
\end{displaymath}
\end{definition}

For $p=2$, this is the notion of $\Lambda^{as}(q)$-set introduced in \cite{LeLiQuRo}. Since the $p$-stable 
norms $[\![ \  ]\!]_p$ dominate the Gaussian norm $[\![ \  ]\!]$, which dominates the norm $\| \ \|_2$, 
every $\Lambda (q)$-set is a $\Lambda^{as}(q)$-set, and every $\Lambda^{as}(q)$-set is a 
$\Lambda^{p-as}(q)$-set for $1 < p < 2$. However:\par
\begin{proposition}\label{Lambda q}
If $\Lambda$ is a $\Lambda (q)$-set, with $ q > 2$, then, for $1 < p < 2$, it is a $\Lambda^{p-as} (r)$-set,  
with $r = \frac{p'}{2}\,q > q$.
\end{proposition}

This is particularly interesting when $\Lambda$ is a ``true'' $\Lambda (q)$-set, i.e. a $\Lambda (q)$-set 
which is not a $\Lambda (s)$-set for any $s > q$ (see \cite{Bourgain3} and \cite{Talagrand2}).\par
\medskip

It is worth to note that this differs from the case $p=2$, since we proved in \cite{LeLiQuRo}, Theorem~4.3, 
that for any $r > 2$, there exist sets which are not $\Lambda^{as} (r)$-sets, though they are 
$\Lambda (q)$-sets for every $q < r$.
\medskip

\noindent{\bf Proof.} Since $\Lambda$ is a $\Lambda (q)$-set, the inverse Fourier transform is continuous 
from $\ell_2 (\Lambda)$ to $L_\Lambda^q$. Since it is trivially continuous from 
$\ell_1 (\Lambda)$ to $L_\Lambda^\infty$, we get, by interpolation, that is is continuous from 
$\ell_p (\Lambda)$ to $L_\Lambda^r$, with $\frac{1}{p} = \frac{1 - \theta}{1} + \frac{\theta}{2}$ and 
$\frac{1}{r} = \frac{1 - \theta}{\infty} + \frac{\theta}{q}$, i.e. $r = \frac{p'}{2}\,q$. Hence, for every 
$f \in F_p (\Lambda)$ ($a \lesssim b$ meaning $a = O\,(b)$ and $a \approx b$ that $a \lesssim b$ and 
$b \lesssim a$):
\begin{align}
\| f \|_r 
& \lesssim \Big( \sum_\Lambda | \widehat f (\gamma) |^p \Big)^{1/p} 
\approx \E\,\Big| \sum_\Lambda Z_\gamma \widehat f (\gamma)\,\gamma \Big| \notag \\
& \leq \E\,\Big\| \sum_\Lambda Z_\gamma \widehat f (\gamma)\,\gamma \Big\|_\infty 
= [\![ f ]\!]_p. 
\tag*{$\square$}
\end{align}

\noindent{\bf Remarks.}  Since Hausdorff-Young inequality asserts that 
$\| f \|_{p'} \leq \| \widehat f \|_p$ and since, as seen above, $\| \widehat f \|_p \lesssim [\![ f ]\!]_p$, 
every subset of $\Gamma$ is a $\Lambda^{p-as} (q)$-set for $q \leq p'$. Hence this notion is only 
interesting for $q > p'$. Note that in Proposition~\ref{Lambda q}, one has $r > p'$.\par
\medskip

For $1 \leq q < p \leq 2$, the same inequality $\| f \|_{q'} \leq \| \widehat f \|_q = \| f \|_{F_q}$ shows that 
every $p$-stable $q$-Rider set is a $\Lambda^{p-as} (q')$-set.
\par\medskip

We will not investigate further this notion here, but only give two results about their thinness.

\begin{proposition}\label{mesh}
For every $\Lambda^{p-as}(q)$-set $\Lambda \subset \Z$, there is $\kappa > 0$ such that, for every 
$N \geq 1$:
\begin{displaymath}
| \Lambda \cap [1, N] | \leq \kappa \, N^{p'/q} \log N.
\end{displaymath}
\end{proposition}
It follows that the set ${\mathbb S}$ of squares is not a $\Lambda^{p-as}(q)$-set of $\Z$ when $q > 2p'$.
\medskip

\noindent{\bf Proof.} It follows the classical one. Write 
$\Lambda \cap [1, N] =\{\lambda_1, \ldots, \lambda_n\}$, and consider the trigonometric polynomial 
$f = \sum_{j=1}^n e_{\lambda_j}$, where $e_{\lambda_j} (t) = \e^{i \lambda_j t}$. By 
Lemma~\ref{majoriz. p-stable norm}, one has $[\![ f ]\!]_p \leq C\,n^{1/p} (\log \lambda_n)^{1/p'}$.\par

Now, on the other hand, $n = f \ast D_N (0)$, where $D_N$ is the $N^{th}$ Dirichlet kernel; hence
\begin{displaymath}
n \leq \| f \|_q  \|D_N\|_{q'} \leq K\,[\![ f ]\!]_p N^{1/q}.
\end{displaymath}

Hence, since $\lambda_n \leq N$, we get $n \leq \kappa\,n^{1/p} (\log N)^{1/p'} N^{1/q}$, and the 
result follows.
\hfill $\square$

\begin{corol}
Let $\alpha$ be an integer $\geq 2$, and $r_\alpha (j)$ is the number of ways to write $j$ as a sum of  
$\alpha$ elements of $\Lambda$. If $\Lambda$ is an $\Lambda^{p-as}(2\alpha)$-set of $\N$, then 
\begin{displaymath}
\frac{1}{n} \sum_{j=1}^n r_\alpha^2 (j) \lesssim n^{\frac{2 - p}{p - 1}} (\log n)^{2\alpha}.
\end{displaymath}
\end{corol}
\noindent{\bf Proof.} We follow Rudin's proof of Theorem~4.5 of \cite{Rudin}. Writing 
$\Lambda = \{n_1, n_2, \ldots\}$, one consider the trigonometric polynomial 
\begin{displaymath}
f (t) =\e^{ i n_1 t} +\cdots + \e^{i n_k t}.
\end{displaymath}
One has
\begin{displaymath}
f^\alpha (t) = r_\alpha(0) + r_\alpha (1) \e^{it} + \cdots\,,
\end{displaymath}
and so $\sum_{j=1}^{n_k} r_\alpha^2(j) \leq \| f \|_{2\alpha}^{2\alpha}$. But:
\begin{displaymath}
\| f \|_{2\alpha} \lesssim [\![ f ]\!]_p \lesssim k^{1/p} (\log n_k)^{1/p'}, 
\end{displaymath} 
and, by Proposition~\ref{mesh}, one has $k \lesssim n_k^{p'/2\alpha} \log n_k$; we get hence 
$[\![ f ]\!]_p^{2\alpha}  \lesssim n_k^{p'/p} (\log n_k)^{2\alpha}$, and the result follows.
\hfill $\square$

\bigskip

Of course other random variables might be used instead of $p$-stable ones. In particular, $1$-stable ones, 
for which one has the quasi-norm (see \cite{MaPi}, page~296):
\begin{displaymath}
[\![ f ]\!]_1 = \sup_{c >0} c\, {\mathbb P}\, \Big( 
\Big\| \sum_\gamma Z_\gamma \widehat f (\gamma)\, \gamma \Big\|_\infty > c \Big),
\end{displaymath}
where $(Z_\gamma)_\gamma$ is an i.i.d. family of $1$-stable random variables. A characterization of 
the continuity of $1$-stable random Fourier series is given in \cite{MaPi2} and \cite{Talagrand}.

\goodbreak


\nobreak

\vbox{\it \footnotesize
\begin{description}
\item[{\rm Pascal Lef\`evre}] \hfill\hbox{}\\ 
Universit\'e d'Artois \\
Laboratoire de Math\'ematiques de Lens EA~2462 \\
F\'ed\'eration CNRS Nord-Pas-de-Calais FR~2956 \\
Facult\'e des Sciences Jean Perrin \\
Rue Jean Souvraz, S.P.\kern 1mm 18 \\ 
62\kern 1mm 307 LENS Cedex
FRANCE \\ 
pascal.lefevre@euler.univ-artois.fr 

\item [{\rm Daniel Li}] \hfill\hbox{}\\ 
Universit\'e d'Artois \\
Laboratoire de Math\'ematiques de Lens EA~2462 \\
F\'ed\'eration CNRS Nord-Pas-de-Calais FR~2956 \\
Facult\'e des Sciences Jean Perrin\\
Rue Jean Souvraz, S.P.\kern 1mm 18\\ 
62\kern 1mm 307 LENS Cedex
FRANCE \\ 
daniel.li@euler.univ-artois.fr

\item [{\rm Herv\'e Queff\'elec}] \hfill\hbox{}\\
Universit\'e des Sciences et Technologies de Lille \\
Labo\-ratoire Paul Painlev\'e U.M.R. CNRS 8524 \\
U.F.R. de Math\'ematiques\\
59\kern 1mm 655 VILLENEUVE D'ASCQ Cedex 
FRANCE \\ 
queff@math.univ-lille1.fr

\item [{\rm Luis Rodr{\'\i}guez-Piazza}] \hfill\hbox{}\\ 
Universidad de Sevilla \\
Facultad de Matem\'aticas, Departamento de An\'alisis Matem\'atico\\ 
Apartado de Correos 1160\\
41\kern 1mm 080 SEVILLA, SPAIN \\ 
piazza@us.es\\
Partially supported by Spanish project MTM2006-05622

\end{description}
}

\end{document}